\documentclass[onecolumn]{autart}
\usepackage{}    

\usepackage{graphicx}           
\usepackage{dcolumn}            
\usepackage{bm}                 

\usepackage{graphics}           
\usepackage{epsfig}             
\usepackage{times}              
\usepackage[fleqn]{amsmath}
\usepackage{amssymb}
\usepackage{extarrows}
\usepackage{amsfonts}
\usepackage{mathrsfs}
\usepackage{fancyhdr}
\usepackage{color}
\usepackage{subfigure}
\usepackage{enumerate}

\allowdisplaybreaks[4] 

\newtheorem{theorem}{Theorem}
\newtheorem{lemma}[theorem]{Lemma}

\newtheorem{definition}{Definition}

\newtheorem{remark}{Remark}

\begin{document}

\begin{frontmatter}

\title{A weak maximum principle-based approach for input-to-state stability analysis of nonlinear parabolic PDEs with boundary disturbances}

\author{Jun Zheng$^{1}$}\ead{zhengjun2014@aliyun.com},
\author{Guchuan Zhu$^{2}$}\ead{guchuan.zhu@polymtl.ca},

\address{$^{1}${School of Mathematics, Southwest Jiaotong University,
        Chengdu 611756, Sichuan, China}\\
        $^{2}$Department of Electrical Engineering, Polytechnique Montr\'{e}al, P.O. Box 6079, Station Centre-Ville, Montreal, QC, Canada H3T 1J4}

\begin{abstract}
In this paper, we introduce a weak maximum principle-based approach to input-to-state stability (ISS) analysis for certain nonlinear partial differential equations (PDEs) with boundary disturbances. Based on the weak maximum principle, a classical result on the maximum estimate of
solutions to linear parabolic PDEs has been extended, which enables the ISS analysis for certain {}{nonlinear} parabolic PDEs with boundary disturbances. To illustrate the application of this method, we establish ISS estimates for a linear reaction-diffusion PDE and a generalized Ginzburg-Landau equation with {}{mixed} boundary disturbances. Compared to some existing methods, the scheme proposed in this paper involves less intensive computations and can be applied to the ISS analysis for a {wide} class of nonlinear PDEs with boundary disturbances.
\end{abstract}
\end{frontmatter}
\section{Introduction}\label{Sec: Introduction}
Originally introduced by Sontag in the late {}{1980's} for finite-dimensional systems, input-to-state stability (ISS) is one of the central notions in the modern theory of nonlinear robust control. It aims at ensuring that disturbances can only induce, in the worst case, a proportional perturbation on the amplitude of the system trajectory. It has been shown that the ISS is an important tool for describing the robust stability property of  infinite dimensional systems governed by partial differential equations (PDEs), and considerable efforts on the study of this subject have been reported in the recent years,  {}{see, e.g., \cite{Argomedo:2012,Argomedo:2013,Dashkovskiy:2010,Dashkovskiy:2013,Dashkovskiy:2013b,Karafyllis:2018iss,Logemann:2013,Mazenc:2011,Mironchenko:2016,Mironchenko:2014b,Mironchenko:2014a,Mironchenko:2018,Prieur:2012,Tanwani:2017,Zheng:201801}.}

It is worth noting that the extension of the notion of ISS to {PDEs} {}{with respect to (w.r.t.)} distributed in-domain disturbances seems to be straightforward, while the investigation on the ISS properties w.r.t. boundary disturbances is much more challenging. The main difficulty lies in the fact that when the disturbances act on the boundaries, the considered systems will usually be governed PDEs with unbounded control operators, which represents an obstacle in the establishment of ISS.

 To tackle this issue, different solutions have been developed recently {}{for ISS analysis of {PDEs} with boundary disturbances, including}:
\begin{enumerate}[(i)]
 \item  {}{the semigroup and admissibility methods for the ISS of \emph{linear parabolic PDEs} \cite{Jacob:2018CDC,Jacob:2019,jacob2016input,Jacob:2018_SIAM,Jacob:2018_JDE}, or \emph{certain semilinear parabolic PDEs} \cite{Schwenninger:2019} inspired by \cite{Zheng:201804}};\label{semigroup}
   \item the approach of spectral decomposition and finite-difference scheme for the ISS of \emph{{}{parabolic PDEs governed by Sturm-Liouville operators}} \cite{Karafyllis:2014,Karafyllis:2016,Karafyllis:2016a,karafyllis2017siam,Karafyllis:2017};\label{spectral decomposition}
        \item {}{the Riesz-spectral approach for ISS of \emph{Riesz-spectral systems} \cite{Lhachemi:201901,Lhachemi:201902};}\label{Riesz-spectral}
   \item the monotonicity-based method for the ISS of \emph{nonlinear {}{parabolic} PDEs with Dirichlet boundary disturbances} \cite{Mironchenko:2019};\label{monotonicity}
   \item the method of De~Giorgi iteration for the ISS w.r.t boundary disturbances of \emph{nonlinear {}{parabolic} PDEs with Dirichlet boundary disturbances} \cite{Zheng:201702,Zheng:201803};\label{iteration}
   \item  variations of Sobolev embedding inequalities for the ISS of \emph{{}{linear or} nonlinear PDEs with certain Robin boundary disturbances} \cite{Zheng:201801,Zheng:201804,Zheng:201803,Zheng:201802};\label{Sobolev inequalities}
   \item  {}{the application of the maximum principle for the ISS of \emph{certain specific  parabolic PDEs with linear boundary conditions} \cite{Zheng:2019a}}.\label{maximum principle}
       \end{enumerate}

{}{The methods in (i) 
  can be used to conduct ISS analysis for certain linear or nonlinear parabolic PDEs (or {other types of} infinite dimensional systems), e.g., parabolic diagonal systems.} 
   {}{The techniques in the category (ii)-(iii) are effective for linear PDEs over $1$-dimensional spatial domains. Whereas, these approaches} may involve heavy computations for PDEs on multidimensional spatial domains or nonlinear PDEs. It should be mentioned that the methods in the categories {}{(iv)-(vi)} can be applied to the establishment of ISS properties for certain nonlinear PDEs with boundary disturbances, or over multidimensional spatial domains. However, none of them can be used to deal with PDEs with boundary conditions as the one given in \eqref{Example 1b} or \eqref{Landau 2b} of this paper. 
{}{Despite} a rapid progress in the study of the ISS for infinite-dimensional systems, \emph{how to establish ISS properties by {}{the} Lyapunov method for nonlinear systems with different boundary disturbances} is still an open research problem.

{}{This paper presents} a weak maximum principle-based approach, {}{which was originally proposed in
\cite{Zheng:2019a}}, for the ISS analysis of {}{a class} of nonlinear parabolic PDEs with different type of  {}{boundary} disturbances. We will show that, combining with the weak maximum principle-based approach, {}{the} Lyapunov method can {}{still} be used to establish ISS estimates for certain nonlinear parabolic PDEs with different type of  {}{boundary} disturbances.  It should be noted that the weak maximum principle was originally used to establish maximum estimates for solutions to elliptic PDEs and has been extensively applied to regularity theory of evolutional PDEs (see, e.g., \cite{Wu2006}). However, the classical method of {}{the} weak maximum principle cannot be directly applied to the ISS analysis for PDE systems. For this reason, we have developed in this work a new estimate for solutions to linear parabolic PDEs (see \eqref{031103} and Remark \ref{Remark 3}). The application of this method is further illustrated in Section \ref{Applications} through the establishment of ISS for linear parabolic PDEs and generalized Ginzburg-Landau equations, with mixed boundary disturbances, respectively. {}{It is noticed that the weak maximum principle was also used in \cite{Mironchenko:2019} to obtain a comparison principle and to establish ISS estimates for monotone parabolic PDEs with Dirichlet boundary disturbances. Different from \cite{Mironchenko:2019}, the weak maximum principle used in this paper and \cite{Zheng:2019a} is for establishing the maximum estimates of the solutions to parabolic PDEs {}{as a crucial step in} ISS analysis for nonlinear parabolic PDEs with different types of boundary disturbances.} An advantage of {}{the} weak maximum principle-based ISS analysis is that it involves much less computations compared to the aforementioned methods. Moreover, it allows dealing with PDEs with more generic types of  {}{boundary} disturbances.
\paragraph{Notation} In this paper, $\mathbb{R}_+$ denotes the set of positive real numbers, $\mathbb{R}_{\geq 0} :=  {}{\{0\}}\cup\mathbb{R}_+$, and $\mathbb{R}_{\leq 0} :=  \mathbb{R}\setminus\mathbb{R}_+$.

Let $\mathcal {K}=\{\gamma : \mathbb{R}_{\geq 0} \rightarrow \mathbb{R}_{\geq 0}|\ \gamma(0)=0,\gamma$ is continuous, strictly increasing$\}$;
$ \mathcal {K}_{\infty}=\{\theta \in \mathcal {K}|\ \lim\limits_{s\rightarrow\infty}\theta(s)=\infty\}$;
$ \mathcal {L}=\{\gamma : \mathbb{R}_{\geq 0}\rightarrow \mathbb{R}_{\geq 0}|\ \gamma$ is continuous, strictly decreasing, $\lim\limits_{s\rightarrow\infty}\gamma(s)=0\}$;
$ \mathcal {K}\mathcal {L}=\{\mu : \mathbb{R}_{\geq 0}\times \mathbb{R}_{\geq 0}\rightarrow \mathbb{R}_{\geq 0}|\ \mu \in \mathcal {K}, \forall t \in \mathbb{R}_{\geq 0}$, and $\mu (s,\cdot)\in \mathcal {L}, \forall s \in {\mathbb{R}_{+}}\}$.

For any $T>0$, let $Q_T=(0,1)\times (0,T)$ and $\partial_pQ_T$ be the parabolic boundary of $Q_T $, i.e., $\partial_pQ_T=(\{0,1\}\times (0,T))\cup ([0,1]\times \{0\})$.

\section{Problem setting and main result}\label{Sec: ISS Analysis}
We consider a class of nonlinear parabolic PDEs with Robin boundary disturbances of the form:
\begin{subequations}\label{0311+1}
\begin{align}
&u_t-a u_{xx}+bu_{x}+cu+h(u)=f(x,t), \ \ (x,t)\in  (0,1)\times \mathbb{R}_+,\\
&\alpha_0u(0,t)-\beta_0u_{x}(0,t)=d_0(t),\ \ t\in \mathbb{R}_+ ,\\
&\alpha_1u(1,t)+\beta_1u_{x}(1,t)=d_1(t),\ \ t\in \mathbb{R}_+,\\
&u(x,0)=\phi(x),\ \ x\in (0,1),
\end{align}
\end{subequations}
where $a\in\mathbb{R}_+, b,c\in \mathbb{R}$, $h$ is a nonlinear term in the equation, $f$ and $d_0,d_1$ are disturbances distributed over the domain $(0,1)$ and on the boundaries $x=0,1$, $\phi$ is the initial value, and $\alpha_i,\beta_i(i=0,1)$ are nonnegative constants.

For the sake of technical development of the scheme proposed in this paper, we always suppose that
\begin{align*}
&f\in {}{C^2}( [0,1]\times \mathbb{R}_{\geq 0} ), d_0,d_1\in {}{C^2}(\mathbb{R}_{\geq 0} ),\phi\in {}{C^2}([0,1 ]),\\
&a\in\mathbb{R}_+, b,c\in \mathbb{R}\ \text{with}\ \frac{b^2}{4a}+c>0,\\
 &\alpha_i,\beta_i\in \mathbb{R}_{\geq 0}\ \text{with}\  \alpha_i+\beta_i>0,i=0,1.
\end{align*}
Furthermore, we assume that
\begin{subequations}\label{condition of a}
\begin{align}
\left\{
\begin{array}{l l}
&-4\left(\dfrac{\alpha_0}{\beta_0}-\dfrac{b}{2a}\right)<\dfrac{b^2}{4a}+c\\
&-\left(\dfrac{\alpha_0}{\beta_0}-\dfrac{b}{2a}\right)\leq a
 \end{array} \right.\ \ \ \ \text{if}\ \ \ \
 \left\{
\begin{array}{l l}
&\beta_0>0\\
&\alpha_0-\dfrac{b}{2a}\beta_0 \leq 0
 \end{array} \right.,\label{+4a}\\
 \left\{
\begin{array}{l l}
&-4\left(\dfrac{\alpha_1}{\beta_1}+\dfrac{b}{2a}\right)<\dfrac{b^2}{4a}+c\\
&-\left(\dfrac{{}{\alpha_1}}{\beta_1}+\dfrac{b}{2a}\right)\leq a
 \end{array} \right.\ \  \ \text{if}\ \ \ \
  \left\{
\begin{array}{l l}
&\beta_1>0\\
&\alpha_1+\dfrac{b}{2a}\beta_1 \leq 0
 \end{array} \right..\label{+4b}
  \end{align}
\end{subequations}
For the nonlinearity, we assume that $h\in {}{C^{2}}( \mathbb{R})$ satisfying
\begin{subequations}\label{condition of h}
\begin{align}
&h(0)=0,\\
& h(|s|)+h(s)\geq 0, \ \forall s\in \mathbb{R},\\
&\frac{b^2}{4a}+c+2h'(s)\geq 0, \ \forall s\in \mathbb{R}_{\leq 0},\\
&h'(s)\geq 0,\ \ \forall s\in \mathbb{R}_+.
\end{align}
\end{subequations}
{}{Besides, we always assume that the following  compatibility conditions are satisfied:
\begin{align}\label{condition 1}
d_{i}'(t)+cd_{i}(t)+\alpha_ih\bigg(\frac{d_{i}(t)}{\alpha_i}\bigg)=\alpha_if(i,t),\forall t\in\mathbb{R}_+,\ \text{for}\ \beta_i=0,i=0\ \text{or}\ 1,
\end{align}
and
{\begin{align}\label{condition 2}
\alpha_i\phi(i)-(-1)^{i}\beta_i\phi_x(i)=d_i(0)=0,\text{for}\ \beta_i\neq0,i=0\ \text{or}\ 1,
\end{align}}}
\begin{remark} {}{We provide some remarks on the structural conditions and well-posedness of a classical solution.}
\begin{enumerate}
  \item[(i)]
If $\frac{|b|}{2a}\leq \frac{\alpha_i}{\beta_i}$ with $\alpha_i>0,\beta_i>0,i=0,1,$ particularly, if $b=0$ with $\alpha_i\geq 0 ,\beta_i>0,i=0,1,$  then conditions in \eqref{condition of a} hold trivially.
\item[(ii)] Examples of $h$ on $\mathbb{R}_{\leq 0} $ include:  (1) $h(s)=\sum\limits_{i=0}^{k}\mu_is^{2i+1}$ with $k$ a positive integer and $\mu_i$ nonnegative constants, $i=1,2,...,2k+1$; (2) $ h(s)=(1+s)\ln(1-s)+s$; (3) $h(s)=\ln (1-ms)-ns$ with constants $ m>0,n>0$ satisfying $m+n\leq \frac{1}{2}\big( \frac{b^2}{4a}+c\big)$.
    \item[(iii)] By \eqref{condition of h}, if we restrict $h$ on $\mathbb{R}_{\geq 0} $, then $h\in \mathcal {K}$.
   \item[(iv)]
   {}{For any $T>0$, \cite[Theorem 6.1, Chapter V]{Ladyzenskaja:1968} (or \cite[Theorem 7.4, Chapter V]{Ladyzenskaja:1968}) guarantees that \eqref{0311+1} with Dirichlet ({or Neumann, or Robin}) boundary conditions admits a unique solution $u\in C^{2,1}(\overline{Q}_T)$. The existence of a unique solution to \eqref{0311+1} with mixed boundary conditions can be {assessed} exactly as in \cite[Chapter~V]{Ladyzenskaja:1968} and obtained by the Leray-Schauder Theorem. 
   It is worth noting that the existence of a global classical solution was also given in \cite{Amann:1989} for quasilinear parabolic equations with boundary conditions {under} a general form.}

\end{enumerate}
\end{remark}
\begin{definition}\label{definition 1}
System~\eqref{0311+1} is said to be input-to-state stable (ISS) in $L^2$-norm w.r.t. boundary disturbances $d_0,d_1$ and in-domain disturbance $f$, if there exist functions $\beta\in \mathcal {K}\mathcal {L}$ and $ \gamma,\gamma_0,\gamma_1\in \mathcal {K}$ such that the solution of \eqref{0311+1} satisfies for any $T>0$:
\begin{align}\label{Eq: ISS def}
    \|u(\cdot,T)\|\leq &\beta( \|{\phi}\|,T)+\gamma_0(\|d_0\|_{L^{\infty}(0,T)})+\gamma_1(\|d_1\|_{L^{\infty}(0,T)})
      +\gamma(\|f\|_{L^{\infty}((0,1)\times (0,T)}).
\end{align}
Moreover, System~\eqref{0311+1} is said to be exponential input-to-state stable (EISS) in $L^2$-norm w.r.t. boundary disturbances $d_0,d_1$ and in-domain disturbance $f$, if there exist {}{positive constants $M_0$ and $\lambda$} such that  $\beta( \|\phi\|,T) \leq{}{ M_0}\|\phi\|e^{-\lambda T}$ in \eqref{Eq: ISS def}.
\end{definition}
{}{The  main result obtained in this paper is the following theorem.}
\begin{theorem}\label{main result}
System \eqref{0311+1} is EISS w.r.t. boundary disturbances $d_0,d_1$ {}{and} in-domain disturbance $f$, having the estimate {}{given in \eqref{Eq: ISS def} with $\beta( \|\phi\|,T)=e^{\frac{|b|}{2a}}\|\phi\|e^{-\lambda T} $,
%
%
$\gamma,\gamma_0,\gamma_1\in  \mathcal {K}$ and $\lambda\in\mathbb{R}_+$ satisfying \eqref{0319+5}.}
\end{theorem}

\begin{remark}
Although the ISS properties obtained in this paper is in the sense of $L^2$-norm, we can treat as in \cite{Zheng:201803} to obtain the ISS in $L^{2p}$-norm for $p\geq 1$. {}{Besides, it is worthy noting that the ISS estimate given in Theorem \ref{main result} have to involve the $L^\infty$-norm of the disturbances due to the usage of the weak maximum principle, which results in the maximum estimates of solutions of parabolic PDEs. }

\end{remark}

\section{Proof of the main result}\label{Sec: Proof}
We introduce first two technical lemmas used in the development of the main result.
\begin{lemma}\label{Lemma 2}\cite{Zheng:201804}
Suppose that $u\in {C^{1}([a,b];\mathbb{R})}$, then we have
 $u^2(c)\leq \dfrac{2}{b-a}\|u\|^2+(b-a)\|u_x\|^2$ \text{ for any } $ c\in[a,b]$.
\end{lemma}
\begin{lemma}\label{maximum principle}(Weak maximum principle \cite[page 237]{Wu2006}) Suppose that $a,b,c,f$ are continuous functions on $Q_T$, and $c\geq 0$ is a bounded function on $Q_T $. If $u\in C^{2,1}(Q_T)\cap C(\overline{Q}_T)$ satisfies
$L_tu:=u_t-au_{xx}+b u_x+cu=f(x,t)\leq 0$ (resp. $\geq 0$) in $  Q_T$, then
\begin{align*}
\max\limits_{\overline{Q}_T}u \leq \max\limits_{\partial_pQ_T}u^+ \ \ \bigg(\text{resp. } \min_{\overline{Q}_T}u \geq -\min_{\partial_pQ_T}u^-\bigg),
\end{align*}
where $u^{+}=\max\{0, u\}$ and $u^{-}=\max\{0, -u\}$.
\end{lemma}

Now we prove the main result of this paper based on the weak maximum principle-approach and {}{the} Lyapunov method.\\
\textbf{Proof of Theorem \ref{main result}} We proceed in four steps.

\textbf{Step 1 (transforming and splitting):} {}{In order to get rid of the term $bu_x$, we use the technique of transforming.} Let $\widetilde{c}= \frac{b^2}{4a}+c$, $\widetilde{\alpha}_0=\alpha_0-\frac{b}{2a}\beta_0 $, $\widetilde{\alpha}_1=\alpha_1+\frac{b}{2a}\beta_1 $ and  $\widetilde{\beta}_i=\beta_i,i=0,1$. Using $\widetilde{u}(x,t)=e^{\frac{-bx}{2a}}u(x,t)$,
 $\widetilde{h}(u)=e^{\frac{-bx}{2a}}h(u)$,
 $\widetilde{\phi}(x)=e^{\frac{-bx}{2a}}\phi(x)$, $\widetilde{f}(x,t)=e^{\frac{-bx}{2a}}f(x,t)$,
 $\widetilde{d}_i(t)=e^{\frac{-bi}{2a}}d_i(t),i=0,1$, we transform \eqref{0311+1} into the following system:
\begin{subequations}\label{Example 2}
\begin{align}
&\widetilde{u}_t-a \widetilde{u}_{xx}+\widetilde{c}\widetilde{u}+\widetilde{h}(e^{\frac{bx}{2a}}(\widetilde{u} ) )=\widetilde{f}(x,t),\ \ (x,t)\in  (0,1)\times \mathbb{R}_+,\\
&\widetilde{\alpha}_0\widetilde{u}(0,t)-\widetilde{\beta}_0\widetilde{u}_{x}(0,t)=\widetilde{d}_0(t),\ \ t\in \mathbb{R}_+ ,\\
&\widetilde{\alpha}_1\widetilde{u}(1,t)+\widetilde{\beta}_1\widetilde{u}_{x}(1,t)=\widetilde{d}_1(t),\ \ t\in \mathbb{R}_+ ,\\
 &\widetilde{u}(x,0)=\widetilde{\phi}(x),\ \ x\in (0,1).
\end{align}
\end{subequations}
{}{In order to deal with the nonlinear term $h(u)$ and apply the weak maximum principle for linear systems, we use the technique of splitting as in \cite{Fabre:1995}, which was applied to ISS analysis of parabolic PDEs based on the approach of De Giorgi iteration in \cite{Zheng:201702}. For} any $T>0$, we split \eqref{Example 2} by $\widetilde{v}+\widetilde{w}=\widetilde{u} $ into two subsystems over the domain $Q_T=(0,1)\times (0,T)$:
\begin{subequations}\label{0311+101}
\begin{align}
&\widetilde{v}_t-a \widetilde{v}_{xx}+\widetilde{c}\widetilde{v}=\widetilde{f}(x,t),\ \ (x,t)\in Q_T, \\
&(\widetilde{\alpha}_0+\widetilde{k}_0)\widetilde{v}(0,t)-\widetilde{\beta}_0\widetilde{v}_{x}(0,t)=\widetilde{d}_0(t),\ \ t\in (0,T), \\
&(\widetilde{\alpha}_1+\widetilde{k}_1)\widetilde{v}(1,t)+\widetilde{\beta}_1\widetilde{v}_{x}(1,t)=\widetilde{d}_1(t),\ \ t\in (0,T),\\
&\widetilde{v}(x,0)=0,\ \ x\in (0,1),
\end{align}
\end{subequations}
and
\begin{subequations}\label{Example 2'.2}
\begin{align}
&\widetilde{w}_t-a \widetilde{w}_{xx}+\widetilde{c}\widetilde{w}+\widetilde{h}(e^{\frac{bx}{2a}}(\widetilde{v}+\widetilde{w} ) )=0,\\
&\widetilde{\alpha}_0\widetilde{w}(0,t)-\widetilde{\beta}_0\widetilde{w}_{x}(0,t)-\widetilde{k}_0\widetilde{v}(0,t)=0,\\
&\widetilde{\alpha}_1\widetilde{w}(1,t)+\widetilde{\beta}_1\widetilde{w}_{x}(1,t)-\widetilde{k}_1\widetilde{v}(1,t)=0,\\
 &\widetilde{w}(x,0)=\widetilde{\phi}(x),\ \  x\in (0,1),
\end{align}
\end{subequations}
 where $\widetilde{k}_i (i=0,1)$ are nonnegative constants satisfying
 \begin{align}\label{$k_i$}
 \left\{
\begin{array}{l l}
&\widetilde{k}_i=0,\ \ \ \ \ \ \ \ \ \ \ \text{if}\ \widetilde{\alpha}_i>0\ \text{or} \ \widetilde{\beta}_i=0\\
&\widetilde{\alpha}_i+\widetilde{k}_i>0,\ \ \text{else}
 \end{array} \right., \ i=0,1.
\end{align}

\textbf{Step 2 (maximum estimate of $\widetilde{v}$):} For \eqref{0311+101}, we establish the following maximum estimate of $\widetilde{v}$ based on {}{the} weak maximum principle (Lemma \ref{maximum principle}):
\begin{align}\label{031103}
    \max\limits_{\overline{Q}_T}|\widetilde{v}|\leq \max\bigg\{ \frac{1}{\widetilde{c}}\sup\limits_{Q_T}|\widetilde{f}|,\frac{1}{\widetilde{\alpha}_0+\widetilde{k}_0} \sup\limits_{ (0,T)} |\widetilde{d}_0|,\frac{1}{\widetilde{\alpha}_1+\widetilde{k}_1} \sup\limits_{ (0,T)} |\widetilde{d}_1|\bigg\}.
\end{align}
Let $\widetilde{M}=\max\bigg\{ \frac{1}{\widetilde{c}}\sup\limits_{Q_T}|\widetilde{f}|,\frac{1}{\widetilde{\alpha}_0+\widetilde{k}_0} \sup\limits_{ (0,T)} |\widetilde{d}_0|,\frac{1}{\widetilde{\alpha}_1+\widetilde{k}_1} \sup\limits_{ (0,T)} |\widetilde{d}_1|\bigg\}$ and $\overline{v}=\widetilde{M}\pm \widetilde{v}$. Then we have
\begin{align*}
 \overline{v}_t-a \overline{v}_{xx}+\widetilde{c} \overline{v}= \widetilde{c}\widetilde{M} \pm (  \widetilde{v}_t-a   \widetilde{v}_{xx}+\widetilde{c} \widetilde{v})=\widetilde{c}{}{\widetilde{M}} \pm \widetilde{f}\geq  \frac{\widetilde{c}}{\widetilde{c}}\sup\limits_{Q_T}|\widetilde{f}|\pm \widetilde{f} \geq 0.
\end{align*}
By {}{the} weak maximum principle (Lemma \ref{maximum principle}), if $\overline{v}$ has a negative minimum, then $\overline{v}$ attains the negative minimum on the parabolic boundary $\partial_pQ_T$. On the other hand, noting that $\overline{v}(x,0)=\widetilde{M}\geq 0$ in $(0,1)$, then $\overline{v}$ attains the negative minimum on $ \{0,1\} \times (0,T)$, i.e., there exists a point $(x_0,t_0)\in \{0,1\} \times (0,T)$, such that $\overline{v}(x_0,t_0) $ is the negative minimum. Thus,
\begin{align*}
&\overline{v}_x(x_0,t_0)\geq 0,\ \ \text{if}\ x_0=0,\\
&\overline{v}_x(x_0,t_0)\leq 0,\ \ \text{if}\ x_0=1.
\end{align*}
Then, at the point $(x_0,t_0)$, we have
\begin{align*}
0>&(\widetilde{\alpha}_0+\widetilde{k}_0)\overline{v}(x_0,t_0)-\widetilde{\beta}_0\overline{v}_{x}(x_0,t_0)
=(\widetilde{\alpha}_0+\widetilde{k}_0)\widetilde{M}\pm \widetilde{d}_0\notag\\
\geq& (\widetilde{\alpha}_0+\widetilde{k}_0)\times \frac{1}{\widetilde{\alpha}_0+\widetilde{k}_0} \sup\limits_{ (0,T)} |\widetilde{d}_0| \pm \widetilde{d}_0
\geq 0,\ \ \text{if}\ x_0=0,
\end{align*}
or
\begin{align*}
0>&(\widetilde{\alpha}_1+\widetilde{k}_1)\overline{v}(x_0,t_0)+\widetilde{\beta}_1\overline{v}_{x}(x_0,t_0)
=(\widetilde{\alpha}_1+\widetilde{k}_1)\widetilde{M}\pm \widetilde{d}_1\notag\\
\geq& (\widetilde{\alpha}_1+\widetilde{k}_1)\times \frac{1}{\widetilde{\alpha}_1+\widetilde{k}_1} \sup\limits_{ (0,T)} |\widetilde{d}_1| \pm \widetilde{d}_1
\geq 0,\ \ \text{if}\ x_0=1,
\end{align*}
both of which lead to a contradiction. Therefore, there must be $\overline{v}\geq0$ in $\overline{Q}_T $, which {}{yields} $|\widetilde{v}|\leq \widetilde{M}$ in $\overline{Q}_T$.

\textbf{Step 3 ($L^2$-estimate of $\widetilde{w}$):} For \eqref{Example 2'.2}, we establish the $L^2$-estimate of $\widetilde{w}$ by {}{the} Lyapunov method and \eqref{031103}, and show that  for any $ T>0$:
 \begin{align}\label{+03185}
    \|\widetilde{w}(\cdot,T)\|\leq& \|\phi\|e^{-\widetilde{\lambda}T}+
\widetilde{\Gamma}\bigg(\sup\limits_{Q_T}|\widetilde{f}| \bigg)+\widetilde{\Gamma}_0\bigg(\sup\limits_{ (0,T)} |\widetilde{d}_0|\bigg)+\widetilde{\Gamma}_1\bigg(\sup\limits_{ (0,T)} |\widetilde{d}_1|\bigg),
\end{align}
where $ \widetilde{\Gamma},\widetilde{\Gamma}_0,\widetilde{\Gamma}_1 \in \mathcal {K}$ are given by
\begin{align*}
\widetilde{\Gamma}(s)=
\widetilde{\mu} s+\widetilde{\sigma} h(\widetilde{\tau} s),
\widetilde{\Gamma}_i(s)=
\widetilde{\mu}_is+ \widetilde{\sigma}_i h(\widetilde{\tau}_i s),i=0,1, \forall s\geq 0,
\end{align*}
with positive constants $\widetilde{\lambda},\widetilde{\mu},\widetilde{\sigma},\widetilde{\tau},\widetilde{\mu}_i,
\widetilde{\sigma}_i,\widetilde{\tau}_i (i=0,1)$ not depending on $T$.

 Indeed, multiplying \eqref{Example 2'.2} with $ \widetilde{w}$ and integrating by parts over $(0,1)$, we obtain
\begin{align}\label{+b}
\frac{1}{2}\frac{\text{d}}{\text{d}t}\|\widetilde{w}\|^2  +a\|\widetilde{w}_x\|^2+\widetilde{c}\|{}{\widetilde{w}}\|^2=a\widetilde{w}_x\widetilde{w}|^{1}_{0}-\int_{0}^1\widetilde{h}(e^{\frac{bx}{2a}}(\widetilde{v}+\widetilde{w} ) )\widetilde{w}\text{d}x.
\end{align}
There are four cases:
 (i) $\widetilde{\beta}_0>0$ and $\widetilde{\beta}_1>0$; (ii) $\widetilde{\beta}_0>0$ and $\widetilde{\beta}_1=0$; (iii) $\widetilde{\beta}_0=0$ and $\widetilde{\beta}_1>0$; (iv) $\widetilde{\beta}_0=0$ and $\widetilde{\beta}_1=0$.

 First, in the case (i), we obtain by \eqref{Example 2'.2}
\begin{align}\label{+b.1}
\widetilde{w}_x\widetilde{w}|^{1}_{0}=&\frac{1}{\widetilde{\beta}_1}\bigg(\widetilde{k}_1\widetilde{v}(1,t)
-\widetilde{\alpha}_1\widetilde{w}(1,t)\bigg)\widetilde{w}(1,t)-\frac{1}{\widetilde{\beta}_0}\bigg(\widetilde{\alpha}_0\widetilde{w}(0,t)
-\widetilde{k}_0\widetilde{v}(0,t)\bigg)\widetilde{w}(0,t)\notag\\
=&\frac{ \widetilde{k}_1}{\widetilde{\beta}_1} \widetilde{v}(1,t)\widetilde{w}(1,t)-\frac{ \widetilde{\alpha}_1}{\widetilde{\beta}_1} \widetilde{w}^2(1,t)-\frac{ \widetilde{\alpha}_0}{\widetilde{\beta}_0} \widetilde{w}^2(0,t)+\frac{ \widetilde{k}_0}{\widetilde{\beta}_0} \widetilde{v}(0,t)\widetilde{w}(0,t)\notag\\
\leq &\frac{ \widetilde{k}_1}{\widetilde{\beta}_1} \bigg(\frac{1}{2\varepsilon_1}\widetilde{v}^2(1,t)+ \frac{\varepsilon_1}{2}\widetilde{w}^2(1,t)\bigg)+\frac{ \widetilde{k}_0}{\widetilde{\beta}_0} \bigg(\frac{1}{2\varepsilon_0}\widetilde{v}^2(0,t)+ \frac{\varepsilon_0}{2}\widetilde{w}^2(0,t)\bigg)-\frac{ \widetilde{\alpha}_1}{\widetilde{\beta}_1} \widetilde{w}^2(1,t)-\frac{ \widetilde{\alpha}_0}{\widetilde{\beta}_0} \widetilde{w}^2(0,t)\notag\\
= & \bigg(\frac{\widetilde{k}_1\varepsilon_1}{2\widetilde{\beta}_1}-
\frac{ \widetilde{\alpha}_1}{\widetilde{\beta}_1}\bigg) \widetilde{w}^2(1,t)+\bigg(\frac{\widetilde{k}_0\varepsilon_0}{2\widetilde{\beta}_0}-
\frac{ \widetilde{\alpha}_0}{\widetilde{\beta}_0}\bigg)\widetilde{w}^2(0,t) + \frac{\widetilde{k}_1}{2\widetilde{\beta}_1\varepsilon_1}\widetilde{v}^2(1,t)
+\frac{\widetilde{k}_0}{2\widetilde{\beta}_0\varepsilon_0}\widetilde{v}^2(0,t),
\end{align}
where $ \varepsilon_0,\varepsilon_1>0$ are small enough.

If $\widetilde{\alpha}_i >0 $ with $i=0\  \text{or}\ 1$, by the choice of $\widetilde{k}_i$ (see \eqref{$k_i$}), it follows {}{that}
\begin{align}\label{+b.1.1}
\bigg(\frac{\widetilde{k}_i\varepsilon_i}{2\widetilde{\beta}_i}-
\frac{ \widetilde{\alpha}_i}{\widetilde{\beta}_i}\bigg)\widetilde{w}^2(i,t)\leq 0 \ \text{and}\ \frac{\widetilde{k}_i}{2\widetilde{\beta}_i\varepsilon_i}\widetilde{v}^2(i,t)=0,\ i=0\  \text{or}\ 1.
\end{align}

 If $\widetilde{\alpha}_i \leq 0 $ with $i=0\  \text{or}\ 1$, by Lemma \ref{Lemma 2}-(i), we have
\begin{align}\label{+b.1.2}
\left(\dfrac{\widetilde{k}_i\varepsilon_i}{2\widetilde{\beta}_i}-
\dfrac{ \widetilde{\alpha}_i}{\widetilde{\beta}_i}\right)\widetilde{w}^2(i,t)\leq \left(\dfrac{\widetilde{k}_i\varepsilon_i}{2\widetilde{\beta}_i}-
\dfrac{ \widetilde{\alpha}_i}{\widetilde{\beta}_i}\right)(2\|\widetilde{w}\|^2+\|\widetilde{w}_x\|^2),\ i=0\  \text{or}\ 1.
\end{align}
Note that as $ \alpha_i,\beta_i\geq 0,i=0,1$, it is impossible that $\widetilde{\alpha}_0\leq 0$ and $\widetilde{\alpha}_1\leq0 $ at the same time. Considering \eqref{+b.1}, \eqref{+b.1.1} and \eqref{+b.1.2}, we have
\begin{align}\label{+b.1'}
\widetilde{w}_x\widetilde{w}|^{1}_{0}
\leq \left\{
\begin{array}{l l}
  0, & \mbox{if $\widetilde{\alpha}_0>0,\widetilde{\alpha}_1>0$, }\\
\left(\dfrac{\widetilde{k}_1\varepsilon_1}{2\widetilde{\beta}_1}-
\dfrac{ \widetilde{\alpha}_1}{\widetilde{\beta}_1}\right)(2\|\widetilde{w}\|^2+\|\widetilde{w}_x\|^2)+ \dfrac{\widetilde{k}_1}{2\widetilde{\beta}_1\varepsilon_1}\widetilde{v}^2(1,t), & \mbox{if  $\widetilde{\alpha}_0>0,\widetilde{\alpha}_1\leq0$,}\\
\left(\dfrac{\widetilde{k}_0\varepsilon_0}{2\widetilde{\beta}_0}-
\dfrac{ \widetilde{\alpha}_0}{\widetilde{\beta}_0}\right)(2\|\widetilde{w}\|^2+\|\widetilde{w}_x\|^2)+ \dfrac{\widetilde{k}_0}{2\widetilde{\beta}_0\varepsilon_0}\widetilde{v}^2(0,t),
 & \mbox{if $\widetilde{\alpha}_0\leq 0,\widetilde{\alpha}_1>0$.}
  \end{array} \right.
\end{align}

 Now we estimate $-\int_{0}^1\widetilde{h}(e^{\frac{bx}{2a}}(\widetilde{v}+\widetilde{w} ) )\widetilde{w}\text{d}x$ in \eqref{+b}. By the Mean Value Theorem, there exists $\xi(x) $ between $ e^{\frac{bx}{2a}}(\widetilde{v}+\widetilde{w} )$ and $e^{\frac{bx}{2a}}\widetilde{v}$ such that
 \begin{align} \label{+b.2}
   -\int_{0}^1\widetilde{h}(e^{\frac{bx}{2a}}(\widetilde{v}+\widetilde{w} ) )\widetilde{w}\text{d}x=& -\int_{0}^1\bigg(\widetilde{h}\big(e^{\frac{bx}{2a}}(\widetilde{v}+\widetilde{w} ) \big)-\widetilde{h}(e^{\frac{bx}{2a}}\widetilde{v} \big)\bigg)\widetilde{w}\text{d}x- \int_{0}^1\widetilde{h}\big(e^{\frac{bx}{2a}}\widetilde{v} \big) \widetilde{w}\text{d}x\notag\\
   =&-\int_{0}^1\widetilde{h}'(\xi)e^{\frac{bx}{2a}} \widetilde{w}\widetilde{w}\text{d}x- \int_{0}^1\widetilde{h}\big(e^{\frac{bx}{2a}}\widetilde{v} \big) \widetilde{w}\text{d}x\notag\\
   =&-\int_{0}^1h'(\xi) \widetilde{w}^2\text{d}x- \int_{0}^1e^{\frac{-bx}{2a}}h\big(e^{\frac{bx}{2a}}\widetilde{v} \big) \widetilde{w}\text{d}x\notag\\
  \leq&-\int_{0}^1h'(\xi) \widetilde{w}^2\text{d}x+ \int_{0}^1e^{\frac{-bx}{2a}}|h\big(e^{\frac{bx}{2a}}\widetilde{v} \big)|\cdot| \widetilde{w}|\text{d}x\notag\\
   \leq & \frac{\widetilde{c}}{2}\|\widetilde{w}\|^2+\int_{0}^1e^{\frac{-bx}{2a}}h\big(|e^{\frac{bx}{2a}}\widetilde{v} |\big)\cdot| \widetilde{w}|\text{d}x\notag\\
         \leq &\frac{\widetilde{c}}{2}\|\widetilde{w}\|^2+\frac{1}{2\varepsilon}e^{\frac{|b|}{a}}\int_{0}^1h^2\big(|e^{\frac{bx}{2a}}\widetilde{v} |\big)\text{d}x+\frac{\varepsilon}{2} \|\widetilde{w}\|^2\notag\\
         \leq &(\frac{\widetilde{c}}{2}+ \frac{\varepsilon}{2})\|\widetilde{w}\|^2+\frac{1}{2\varepsilon}e^{\frac{|b|}{a}}h^2\bigg(e^{\frac{|b|}{2a}}\max_{[0,1]}|\widetilde{v} |\bigg),
\end{align}
where we used \eqref{condition of h}, and $ \varepsilon >0$ will be chosen later.

 We deduce {from} \eqref{+b}, \eqref{+b.1'} and \eqref{+b.2} {}{that}
\begin{align}\label{+03181}
\frac{1}{2}\frac{\text{d}}{\text{d}t}\|\widetilde{w}\|^2  +a\|\widetilde{w}_x\|^2+\widetilde{c}\|{}{\widetilde{w}}\|^2
\leq \left\{
\begin{array}{l l}
  \left(\dfrac{\widetilde{c}}{2}+ \dfrac{\varepsilon}{2}\right)\|\widetilde{w}\|^2+\widetilde{C}_h(t), &  \ \mbox{if $\widetilde{\alpha}_0>0,\widetilde{\alpha}_1>0$, }\\
 \left(2\widetilde{C}_1+\dfrac{\widetilde{c}}{2}+ \dfrac{\varepsilon}{2}\right)\|\widetilde{w}\|^2+\widetilde{C}_1\|\widetilde{w}_x\|^2
 +\widetilde{C}_h(t)+\widetilde{V}_1(t), & \  \mbox{ if $\widetilde{\alpha}_0>0,\widetilde{\alpha}_1\leq0$,}\\
 \left(2\widetilde{C}_0+\dfrac{\widetilde{c}}{2}+ \dfrac{\varepsilon}{2}\right)\|\widetilde{w}\|^2+\widetilde{C}_0\|\widetilde{w}_x\|^2
 +\widetilde{C}_h(t)+\widetilde{V}_0(t),
 &  \ \mbox{if $\widetilde{\alpha}_0\leq 0,\widetilde{\alpha}_1>0$,}
  \end{array} \right.
\end{align}
 where
    \begin{align*}
 & \widetilde{C}_i=\frac{\widetilde{k}_i\varepsilon_i}{2\widetilde{\beta}_i}-
\frac{ \widetilde{\alpha}_i}{\widetilde{\beta}_i},\widetilde{V}_i(t)=\frac{\widetilde{k}_i}{2\widetilde{\beta}_i\varepsilon_i}\widetilde{v}^2(i,t),i=0,1,\widetilde{C}_h(t)=\frac{1}{2\varepsilon}e^{\frac{|b|}{a}}
 h^2\bigg(e^{\frac{|b|}{2a}}\max\limits_{x\in[0,1]}|\widetilde{v}(x,t) |\bigg).
\end{align*}

 Note that by \eqref{condition of a}, one can choose $ \varepsilon,\varepsilon_i,i=0,1,$ small enough, such that
 \begin{subequations}\label{0319+2}
  \begin{align}
 & \frac{\widetilde{c}}{2}+\frac{\varepsilon}{2}<\widetilde{c},\\
&2\widetilde{C}_1+\frac{\widetilde{c}}{2}+ \frac{\varepsilon}{2}< \widetilde{c},\ \text{and}\ \widetilde{C}_1\leq a, \ \text{if}\ \widetilde{\alpha}_0>0,\widetilde{\alpha}_1\leq0,\\
&2\widetilde{C}_0+\frac{\widetilde{c}}{2}+ \frac{\varepsilon}{2}< \widetilde{c},\ \text{and}\ \widetilde{C}_1\leq a, \ \text{if}\ \widetilde{\alpha}_0\leq 0,\widetilde{\alpha}_1>0.
\end{align}
\end{subequations}
 Then we infer from \eqref{+03181} {}{that}
 \begin{align}\label{+03182}
\frac{\text{d}}{\text{d}t}\|\widetilde{w}(\cdot,t)\|^2
      \leq& -2\widetilde{\lambda} \|\widetilde{w}(\cdot,t)\|^2+2( \widetilde{C}_h(t)+\widetilde{V}_1(t)+\widetilde{V}_0(t)),
      \ \forall t\in (0,T),
\end{align}
 where $ \widetilde{\lambda}=\min\limits_{i=0,1}\{\widetilde{\lambda}_i  \}, \widetilde{\lambda}_i =\widetilde{c}-\big(2\widetilde{C}_i+\frac{\widetilde{c}}{2}+ \frac{\varepsilon}{2} \big)  >0, i=0,1$.

By Gronwall's inequality, it follows {}{that}
\begin{align}\label{+1}
\|\widetilde{w}(\cdot,t)\|^2 \leq & \|\widetilde{\phi}\|^2e^{-2 \widetilde{\lambda} t}+2\int_{0}^t
( \widetilde{C}_h(s)+\widetilde{V}_1(s)+\widetilde{V}_0(s))
e^{-2\widetilde{\lambda}(t-s)}\text{d}s
\notag\\
\leq & \|\phi\|^2e^{-2\widetilde{\lambda}t}+2\max\limits_{0\leq s\leq t}( \widetilde{C}_h(s)+\widetilde{V}_1(s)+\widetilde{V}_0(s))\int_{0}^te^{-2\widetilde{\lambda}(t-s)}\text{d}s
\notag\\
\leq& \|\phi\|^2e^{-2\widetilde{\lambda}t}+\frac{2 }{\widetilde{\lambda}}\max\limits_{0\leq s\leq t}(\widetilde{C}_h(s)+\widetilde{V}_1(s)+\widetilde{V}_0(s)).
\end{align}
By \eqref{031103}, it follows that
\begin{align*}
    \max\limits_{x\in[0,1]}|\widetilde{v}(x,t)|\leq \max\bigg\{ \frac{1}{\widetilde{c}}\sup\limits_{Q_T}|\widetilde{f}|,\frac{1}{\widetilde{\alpha}_0+\widetilde{k}_0} \sup\limits_{ (0,T)} |\widetilde{d}_0|,\frac{1}{\widetilde{\alpha}_1+\widetilde{k}_1} \sup\limits_{ (0,T)} |\widetilde{d}_1|\bigg\},
    \end{align*}
which yields
\begin{align}\label{+03183}
\widetilde{V}_i(s)\leq \frac{\widetilde{k}_i}{2\widetilde{\beta}_i\varepsilon_i}\bigg( \frac{1}{\widetilde{c}^2}\sup\limits_{Q_T}|\widetilde{f}|^2+\frac{1}{(\widetilde{\alpha}_0+\widetilde{k}_0)^2} \sup\limits_{ (0,T)} |\widetilde{d}_0|^2+\frac{1}{(\widetilde{\alpha}_1+\widetilde{k}_1)^2} \sup\limits_{ (0,T)} |\widetilde{d}_1|^2\bigg),
\end{align}
and
\begin{align}\label{+03184}
\widetilde{C}_h(s)=&\frac{1}{2\varepsilon}e^{\frac{|b|}{a}}
 h^2\bigg(e^{\frac{|b|}{2a}}\max\limits_{x\in[0,1]}|\widetilde{v}(x,s) |\bigg)\notag\\
 \leq& \frac{1}{2\varepsilon}e^{\frac{|b|}{a}} h^2\bigg(\max\bigg\{e^{\frac{|b|}{2a}}\big( \frac{1}{\widetilde{c}}\sup\limits_{Q_T}|\widetilde{f}|,
 \frac{1}{\widetilde{\alpha}_0+\widetilde{k}_0} \sup\limits_{ (0,T)} |\widetilde{d}_0|,\frac{1}{\widetilde{\alpha}_1+\widetilde{k}_1} \sup\limits_{ (0,T)} |\widetilde{d}_1|\big)\bigg\}\bigg)\notag\\
 \leq &\frac{1}{2\varepsilon}e^{\frac{|b|}{a}} \bigg( h^2\bigg(e^{\frac{|b|}{2a}}\frac{1}{\widetilde{c}}\sup\limits_{Q_T}|\widetilde{f}|\bigg)+
 h^2\bigg(e^{\frac{|b|}{2a}}\frac{1}{\widetilde{\alpha}_0+\widetilde{k}_0} \sup\limits_{ (0,T)} |\widetilde{d}_0|\bigg)+h^2\bigg(e^{\frac{|b|}{2a}} \frac{1}{\widetilde{\alpha}_1+\widetilde{k}_1} \sup\limits_{ (0,T)} |\widetilde{d}_1|\bigg)\bigg).
\end{align}
By \eqref{+1}, \eqref{+03183} and \eqref{+03184}, we find that for any $t\in (0,T)$, it holds {}{that}
\begin{align*}
\|\widetilde{w}(\cdot,t)\|\leq& \|\phi\|e^{-\widetilde{\lambda}t}+
\widetilde{\Gamma}\left(\sup\limits_{Q_T}|\widetilde{f}| \right)+\widetilde{\Gamma}_0\left(\sup\limits_{ (0,T)} |\widetilde{d}_0|\right)+\widetilde{\Gamma}_1\left(\sup\limits_{ (0,T)} |\widetilde{d}_1|\right),
%
%
%
%
\end{align*}
where $ \widetilde{\Gamma},\widetilde{\Gamma}_0,\widetilde{\Gamma}_1 \in \mathcal {K}$ are given by
\begin{align*}
\widetilde{\Gamma}(s)=&\frac{1}{\widetilde{c}}
\left(\sqrt{\frac{\widetilde{k}_0}{\widetilde{\lambda}
\widetilde{\beta}_0\varepsilon_0}}+
\sqrt{\frac{\widetilde{k}_1}{\widetilde{\lambda}
\widetilde{\beta}_1\varepsilon_1}}\right) s+\sqrt{\frac{1 }{\widetilde{\lambda} \varepsilon}}e^{\frac{|b|}{2a}}h\left(\frac{1}{\widetilde{c}}e^{\frac{|b|}{2a}} s\right),\ \forall s\geq 0,\\
\widetilde{\Gamma}_0(s)=&\frac{1}{\widetilde{\alpha}_0+\widetilde{k}_0}\left(\sqrt{\frac{\widetilde{k}_0}{\widetilde{\lambda}
\widetilde{\beta}_0\varepsilon_0}}+
\sqrt{\frac{\widetilde{k}_1}{\widetilde{\lambda}
\widetilde{\beta}_1\varepsilon_1}}\right)s+
\sqrt{\frac{1 }{\widetilde{\lambda} \varepsilon}}e^{\frac{|b|}{2a}}h\left(\frac{1}{\widetilde{\alpha}_0+\widetilde{k}_0}
e^{\frac{|b|}{2a}} s\right),\ \forall s\geq 0,\\
\widetilde{\Gamma}_1(s)=&\frac{1}{\widetilde{\alpha}_1+\widetilde{k}_1}\left(\sqrt{\frac{\widetilde{k}_0}{\widetilde{\lambda}
\widetilde{\beta}_0\varepsilon_0}}+
\sqrt{\frac{\widetilde{k}_1}{\widetilde{\lambda}
\widetilde{\beta}_1\varepsilon_1}}\right)s+
\sqrt{\frac{1 }{\widetilde{\lambda} \varepsilon}}e^{\frac{|b|}{2a}}h\left(\frac{1}{\widetilde{\alpha}_1+\widetilde{k}_1}e^{\frac{|b|}{2a}} s\right),\ \forall s\geq 0.
\end{align*}
Finally, by the continuity of $\widetilde{w}(x,t)$ and $e^{-\widetilde{\lambda}t}$ at $t=T$, it follows {}{that}
\begin{align}\label{+03191}
    \|\widetilde{w}(\cdot,T)\|\leq& \|\phi\|e^{-\widetilde{\lambda}T}+
\widetilde{\Gamma}\left(\sup\limits_{Q_T}|\widetilde{f}| \right)+\widetilde{\Gamma}_0\left(\sup\limits_{ (0,T)} |\widetilde{d}_0|\right)+\widetilde{\Gamma}_1\left(\sup\limits_{ (0,T)} |\widetilde{d}_1|\right),
\end{align}
which completes the proof of Case~(i).

In the case of (ii) {}{(or (iii))}, i.e., $\widetilde{\beta}_i>0$ and $\widetilde{\beta}_{1-i}=0$, $i=0$ {}{(or $1$)}, it suffices to set $ \widetilde{k}_{1-i}=\frac{\widetilde{k}_{1-i}}{\widetilde{\beta}_{1-i}}=\frac{\widetilde{\alpha}_{1-i}}{\widetilde{\beta}_{1-i}}=0$, $i=0$ {}{(or $1$),} in the proof of Case~(i), and to obtain \eqref{+03191}
 with $ \widetilde{\Gamma},\widetilde{\Gamma}_0,\widetilde{\Gamma}_1 \in \mathcal {K}$ given by
 \begin{subequations}\label{0320+1}
\begin{align}
\widetilde{\Gamma}(s)=&\frac{1}{\widetilde{c}}
\sqrt{\frac{\widetilde{k}_i}{\widetilde{\lambda}
\widetilde{\beta}_i\varepsilon_i}}s+\sqrt{\frac{1 }{\widetilde{\lambda} \varepsilon}}e^{\frac{|b|}{2a}}h\bigg(\frac{1}{\widetilde{c}}e^{\frac{|b|}{2a}} s\bigg),\ \forall s\geq 0,\\
\widetilde{\Gamma}_0(s)=&\frac{1}{\widetilde{\alpha}_i+\widetilde{k}_i}\sqrt{\frac{\widetilde{k}_i}{\widetilde{\lambda}
\widetilde{\beta}_i\varepsilon_i}}
s+
\sqrt{\frac{1 }{\widetilde{\lambda} \varepsilon}}e^{\frac{|b|}{2a}}h\bigg(\frac{1}{\widetilde{\alpha}_i+\widetilde{k}_i}e^{\frac{|b|}{2a}} s\bigg),\ \forall s\geq 0,\\
\widetilde{\Gamma}_1(s)=&\frac{1}{\widetilde{\alpha}_{1-i}}\sqrt{\frac{\widetilde{k}_i}{\widetilde{\lambda}
\widetilde{\beta}_i\varepsilon_i}}s+
\sqrt{\frac{1 }{\widetilde{\lambda} \varepsilon}}e^{\frac{|b|}{2a}}h\bigg(\frac{1}{\widetilde{\alpha}_{1-i}}e^{\frac{|b|}{2a}} s\bigg),\ \forall s\geq 0,
\end{align}
\end{subequations}
where $i=0$ {}{(or $1$)}.

Similarly, in Case~(iv), i.e., $\widetilde{\beta}_0=\widetilde{\beta}_{1}=0$, we obtain \eqref{+03191}
 with $ \widetilde{\Gamma},\widetilde{\Gamma}_0,\widetilde{\Gamma}_1 \in \mathcal {K}$ given by
  \begin{subequations}\label{0320+2}
\begin{align}
\widetilde{\Gamma}(s)=&\sqrt{\frac{1 }{\widetilde{\lambda} \varepsilon}}e^{\frac{|b|}{2a}}h\bigg(\frac{1}{\widetilde{c}}e^{\frac{|b|}{2a}} s\bigg),\ \forall s\geq 0,\\
\widetilde{\Gamma}_0(s)=&
\sqrt{\frac{1 }{\widetilde{\lambda} \varepsilon}}e^{\frac{|b|}{2a}}h\bigg(\frac{1}{\widetilde{\alpha}_0}e^{\frac{|b|}{2a}} s\bigg),\ \forall s\geq 0,\\
\widetilde{\Gamma}_1(s)=&
\sqrt{\frac{1 }{\widetilde{\lambda} \varepsilon}}e^{\frac{|b|}{2a}}h\bigg(\frac{1}{\widetilde{\alpha}_{1}}e^{\frac{|b|}{2a}} s\bigg),\ \forall s\geq 0.
\end{align}
 \end{subequations}

\textbf{Step 4:}
We establish the ISS estimate in $L^2$-norm for \eqref{0311+1}. Indeed, noting that $u=e^{\frac{bx}{2a}}\widetilde{u}(x,t)$ and $\widetilde{u}=\widetilde{v}+\widetilde{w} $, by \eqref{031103} and \eqref{+03185}, we get for any $T>0$
\begin{align*}
\|u(\cdot,T)\|\leq &e^{\frac{|b|}{2a}}\|\widetilde{u}(\cdot,T)\|\\
\leq &e^{\frac{|b|}{2a}}(\|\widetilde{v}(\cdot,T)\| +\|\widetilde{w}(\cdot,T)\|)\\
\leq &e^{\frac{|b|}{2a}}\|\phi\|e^{-\widetilde{\lambda}T}+ e^{\frac{|b|}{2a}}\left( \frac{1}{\widetilde{c}}\sup\limits_{Q_T}|\widetilde{f}|+ \widetilde{\Gamma}\left(\sup\limits_{Q_T}|\widetilde{f}| \right)\right)+e^{\frac{|b|}{2a}}\left( \frac{1}{\widetilde{\alpha}_0+\widetilde{k}_0} \sup\limits_{ (0,T)} |\widetilde{d}_0|+  \widetilde{\Gamma}_0\left(\sup\limits_{ (0,T)} |\widetilde{d}_0|\right) \right)\\
&+e^{\frac{|b|}{2a}}\left(   \frac{1}{\widetilde{\alpha}_1+\widetilde{k}_1} \sup\limits_{ (0,T)} |\widetilde{d}_1|+\widetilde{\Gamma}_1\left(\sup\limits_{ (0,T)} |\widetilde{d}_1|\right)\right)\\
\leq &e^{\frac{|b|}{2a}}\|\phi\|e^{-{\lambda}T}+ \gamma\left( \sup\limits_{Q_T}|f|\right)+  {\gamma}_0\left(\sup\limits_{ (0,T)} |d_0|\right) +{\gamma}_1\left(\sup\limits_{ (0,T)} |d_1|\right),\end{align*}
where
\begin{subequations}\label{0319+5}
\begin{align}
&\gamma(s)= e^{\frac{|b|}{2a}}\left( \frac{1}{\widetilde{c}}e^{\frac{|b|}{2a}}s+ \widetilde{\Gamma}\left(e^{\frac{|b|}{2a}}s \right)\right),\ \forall s\geq 0,\\
&{\gamma}_0(s)=e^{\frac{|b|}{2a}}\left( \frac{1}{\widetilde{\alpha}_0+\widetilde{k}_0}s+  \widetilde{\Gamma}_0(s) \right),\ \forall s\geq 0,\\
&{\gamma}_1(s)=e^{\frac{|b|}{2a}}\left(   \frac{1}{\widetilde{\alpha}_1+\widetilde{k}_1} e^{\frac{|b|}{2a}}s+\widetilde{\Gamma}_1\left(e^{\frac{|b|}{2a}}s\right)\right),\ \forall s\geq 0,\\
&\lambda=\widetilde{\lambda}=\min\limits_{i=0,1}\left\{\widetilde{\lambda}_i  \right\}=\min\limits_{i=0,1}\left\{ \widetilde{c}-\left(2\widetilde{C}_i+\frac{\widetilde{c}}{2}+ \frac{\varepsilon}{2} \right) \right\}  >0, i=0,1,\label{0319+5d}
\end{align}
\end{subequations}
with
\begin{align}\label{0320+3}
\widetilde{C}_i=\left\{
\begin{array}{l l}
    \dfrac{\widetilde{k}_i\varepsilon_i}{2\widetilde{\beta}_i}-
\dfrac{ \widetilde{\alpha}_i}{\widetilde{\beta}_i}, &  \mbox{\ if $\widetilde{\beta}_i>0$}\\
    0,&  \mbox{\ if $\widetilde{\beta}_i=0$}
      \end{array} \right., \ i=0, 1,
      \end{align}
$\widetilde{k}_i$ determined by \eqref{$k_i$}, and $ \varepsilon,\varepsilon_i$ determined by \eqref{0319+2}, $i=0, 1.$ \hfill $\blacksquare$

\begin{remark}
{}{It is worth noting that the requirement of $C^2$-continuity on $h,f,d_0,d_1,\phi$ {and the compatibility conditions \eqref{condition 1} and \eqref{condition 2} are} only for ensuring the existence of classical solutions of the system \eqref{0311+1} and the subsystem \eqref{0311+101}, {which can eventually} be relaxed for the ISS analysis if weak solutions of \eqref{0311+1} are considered. Indeed, it suffices to impose certain conditions to guarantee a weak solution of \eqref{0311+1} and a classical solution of the linear $\widetilde{v}$-subsystem \eqref{0311+101} in Section 3.     
For example, we can {weaken} the assumptions on $f,d_0,d_1$ to be ``$f\in C^{l,\frac{l}{2}}( [0,1]\times \mathbb{R}_{\geq 0} ), d_0,d_1\in C^{\frac{l}{2}+1}(\mathbb{R}_{\geq 0} )$ with some constant $l>0$'', {and relax the compatibility conditions \eqref{condition 1} and \eqref{condition 2} to be:
\begin{align*}
d_{i}'(t)+cd_{i}(t)&=\alpha_if(i,t),\forall t\in\mathbb{R}_+,\ \text{for}\ \beta_i=0,i=0\ \text{or}\ 1,\\
d_i(0)&=0,\text{for}\ \beta_i\neq0,i=0\ \text{or}\ 1,
\end{align*}} then \eqref{0311+101} has a unique solution belonging to $C^{l+2,\frac{l}{2}+1}( \overline{Q}_T ) $, see, e.g., \cite[Theorem 5.2 and 5.3, Chapter IV]{Ladyzenskaja:1968}. For a  weak solution of \eqref{Example 2'.2}, we can relax the assumptions on $h,\phi$ to be ``$h\in C^1(\mathbb{R})$ and $\phi\in L^2(0,1)$''. Noting that by the structural condition \eqref{condition of h}, we always have
\begin{align*}
    -h(s)s=-h'(\xi)s\leq \frac{1}{2}\bigg(\frac{b^2}{4a}+c\bigg)|s|, \forall s\in \mathbb{R},
\end{align*}
where $\xi$ is between $0$ and $s$. Then the existence of a unique weak solution of  \eqref{0311+1} can be {obtained by proceeding} exactly as in \cite[Theorem 6.5, 6.39 and 6.46]{Lieberman:2005} {with} the usual approximation  argument based on \textit{a priori} estimates.}
\end{remark}

\begin{remark}\label{Remark 3}
{}{A crucial step in the proof lies in the maximum estimates for the solutions of parabolic equations. It should be mentioned that the} result of the maximum estimate given by \eqref{031103} is {}{different from the classical maximum estimate of solutions to parabolic equations. For example}, in \cite[page 239]{Wu2006}, a classical maximum estimate of the solutions to linear parabolic equations {}{in a finite time interval $(0, T)$} is given as below:
\begin{itemize}
  \item[] {\emph{{}{Assume that $c\geq 0$ is bounded in $Q_T $. If $u\in C^{2,1}(Q_T)\cap C(\overline{Q}_T)$ satisfies} $L_tu:=u_t-a u_{xx}+bu_{x}+cu=f$ in $Q_T $, then
$
\max\limits_{\overline{Q}_T}|u|\leq \sup\limits_{\partial_pQ_T}|u|+T\sup\limits_{Q_T}|f|.
$}}
\end{itemize}
A main improvement obtained in \eqref{031103} is that the coefficients of {}{$\sup\limits_{Q_T}|\widetilde{f}|,\sup\limits_{Q_T}|\widetilde{d}_0|$ and $\sup\limits_{Q_T}|\widetilde{d}_1|$} do not depend on $T$. It is an essential feature for the establishment of ISS properties for PDEs with boundary disturbances.
\end{remark}

\begin{remark}
It should be mentioned that Lemma \ref{maximum principle} cannot be used directly to establish the maximum estimate for the solution of \eqref{0311+101} if $\widetilde{\alpha}_i\leq 0,i=0$ or $1$. To overcome this difficulty, additional terms $\widetilde{k}_i\widetilde{v}(i,t) (i=0,1)$ are added on the boundaries to guarantee that $\widetilde{\alpha}_i+ \widetilde{k}_i\widetilde{v}>0$ when we use the technique of splitting. Thus, $\widetilde{k}_i\widetilde{v}(i,t)$ can be seen as a stabilizing feedback control with boundary disturbances or nonhomogeneous boundary conditions. The idea of using this type of compensation comes from the so-called Penalty Method in mathematics (see, e.g., \cite{Kinderlehrer:1980}) and its applications in singular free boundary problems (see, Section~2 and Section~3 in \cite{Zhengmath:2017}), and the ISS for nonlinear PDEs with nonlinear boundary disturbances (see \cite{Zheng:201802}).
\end{remark}

\begin{remark}  {}{If the boundary disturbance $\widetilde{d}_1(t)$ is replaced by $\int_{0}^1 k(1,y)\widetilde{v}(y,t)\text{d}y+\widetilde{d}_1(t)$, where $k\in C^2([0,1]\times [0,1])$,} then the condition $\frac{b^2}{4a}+c>0$ can be weakened in the proof of Theorem \ref{main result}. Indeed, using the technique of backstepping  {}{and Volterra integral transformation (see, e.g., \cite{Liu:2003,Liu:2000,Smyshlyaev:2004}), one can transform \eqref{0311+101} into a new system of the following form:}
\begin{align*}
&\widehat{v}_t-a \widehat{v}_{xx}+\widehat{c} \widehat{v}=\widetilde{f}(x,t),\\
&(\widetilde{\alpha}_0+\widetilde{k}_0)\widehat{v}(0,t)-\widetilde{\beta}_0\widehat{v}_{x}(0,t)=\widehat{d}_0(t),\\
&(\widetilde{\alpha}_1+\widetilde{k}_1)\widehat{v}(1,t)+\widetilde{\beta}_1\widehat{v}_{x}(1,t)=\widehat{d}_1(t),\\
 &\widehat{v}(x,0)=0,\ \ x\in(0,1),
\end{align*}
where $\widehat{c}>0$ can be an arbitrary constant. Then using the techniques in this paper and arguing as \cite[Section V]{Zheng:201702}, one may establish the maximum estimate for the solution of the above system, and hence for \eqref{0311+101}.
\end{remark}
\begin{remark}
It is of interest to study ISS properties based on the weak maximum principle for PDEs with a general nonlinearity $h(u)$ distributed over the domain and nonlinearities on the boundaries, e.g.:
\begin{align*}
&L_tu:=u_t-a u_{xx}+bu_{x}+cu=h(u),\ \ (x,t)\in (0,1)\times \mathbb{R}_+,\\
 &u_{x}(0,t)-\Psi_0 (u(0,t))=d_0(t),\ \ t\in \mathbb{R}_+,\\
&u_{x}(1,t)+\Psi_1 (u(1,t))=d_1(t),\ \  t\in \mathbb{R}_+,\\
&u(x,0)=\phi(x),\ \ x\in (0,1),
 \end{align*}
where $h,\Psi_0,\Psi_1$ are nonlinear functions.

{}{Moreover, due} to the nonlinearities on the boundaries, much more arguments are needed to obtain the maximum estimates for solutions of the above equation when the weak maximum principle is used. {}{Furthermore, as shown in  \cite{Zheng:2020}, the method proposed in this paper can be applied to a wider class of nonlinear higher dimensional PDEs with variable coefficients and nonlinear boundary conditions. Finally, one may find that the approach {presented in this work can be extended to ISS analysis of other} nonlinear (abstract) systems by estimating the solution $\widetilde{w}$ of the subsystem \eqref{Example 2'.2} in an abstract form as in \cite[Section 3]{Schwenninger:2019}. }
\end{remark}


\section{Illustrative Examples}\label{Applications}
\subsection{1-$D$ linear reaction-diffusion equation}
\vspace*{-6pt}
We consider the following 1-$D$ linear reaction-diffusion PDE with mixed boundary conditions:
\begin{subequations}\label{Example 1}
\begin{align}
&u_t-a u_{xx}+bu_{x}+cu=f(x,t),\ \ (x,t)\in  (0,1)\times \mathbb{R}_+,\\
&u(0,t)=d_0(t),u_x(1,t)=-K_1u(1,t)+d_1(t),\ \ t\in \mathbb{R}_+,\label{Example 1b}\\
  &u(x,0)=\phi(x),\ \ x\in (0,1).
\end{align}
\end{subequations}
 Assume that {}{
$a,K_1\in\mathbb{R}_+,b,c\in \mathbb{R}\ \text{with}\ \frac{b^2}{4a}+c>0,K_1>\frac{|b|}{2a}$.}
  For \eqref{Example 1}, we have
 \begin{align*}
&\widetilde{c}=\frac{b^2}{4a}+c,h(s)\equiv0,\forall s\in \mathbb{R},\\
&\alpha_0=1,\beta_0=0,\alpha_1=K_1,\beta_1=1,\\
&\widetilde{\alpha}_0=\alpha_0-\frac{b}{2a}\beta_0=1,\widetilde{\beta}_0=\beta_0=0, \\ &\widetilde{\alpha}_1=\alpha_1+\frac{b}{2a}\beta_1=K_1+\frac{b}{2a}>0,\widetilde{\beta}_1=\beta_1=1.
\end{align*}
Then \eqref{+4b} and \eqref{condition of h} hold. Therefore, \eqref{Example 1} is EISS.

Furthermore, according to \eqref{0320+3} and \eqref{0319+5d} (and noting that $\widetilde{k}_0=\widetilde{k}_1=0$ in \eqref{$k_i$}), it follows {}{that}
\begin{align*}
&\widetilde{C}_0=0,\widetilde{C}_1=-\frac{\widetilde{\alpha}_1}{\widetilde{\beta}_1}=-K_1-\frac{b}{2a}<0,\\
&\lambda=\widetilde{\lambda}=\min\left\{ \frac{\widetilde{c}}{2}- \frac{\varepsilon}{2} , \widetilde{c}-\left(2\widetilde{C}_1+\frac{\widetilde{c}}{2}+ \frac{\varepsilon}{2} \right) \right\} 
= \frac{1}{2}\left(\frac{b^2}{4a}+c \right)+2\left(K_1+\frac{b}{2a}\right)- \frac{\varepsilon}{2}>0,
\end{align*}
where we choose $ \varepsilon>0$ such that $ 2(K_1+\frac{b}{2a})- \frac{\varepsilon}{2}>0$.

By \eqref{0320+1} and $\widetilde{k}_1=0$, we have
\begin{align*}
\widetilde{\Gamma}(s)=\widetilde{\Gamma}_0(s)=\widetilde{\Gamma}_1(s)=0,\ \forall s\geq 0,
\end{align*}
which implies {}{that} in \eqref{0319+5}
\begin{align*}
\gamma(s)= \frac{4ae^{\frac{|b|}{a}}}{b^2+4ac}s,{\gamma}_0(s)=e^{\frac{|b|}{2a}}s, \ {\gamma}_1(s)= \frac{2ae^{\frac{|b|}{2a}} }{2aK_1+b}s,\ \forall s\geq 0.
\end{align*}
Finally, the ISS estimate of \eqref{Example 1} is given by
\begin{align*}
\|u(\cdot,T)\|\leq  &e^{\frac{|b|}{2a}}\|\phi\|e^{-\big(\frac{1}{2}(\frac{b^2}{4a}+c )+2(K_1+\frac{b}{2a})- \frac{\varepsilon}{2}\big)T}\notag\\
&+ \frac{4ae^{\frac{|b|}{a}}}{b^2+4ac} \sup\limits_{(0,1)\times(0,T)}|f|+ e^{\frac{|b|}{2a}}\sup\limits_{ (0,T)} |d_0| +\frac{2ae^{\frac{|b|}{2a}} }{2aK_1+b}\sup\limits_{ (0,T)} |d_1|,\ \ \forall T>0.
\end{align*}

\subsection{Ginzburg-Landau equations with real coefficients}
Consider the generalized Ginzburg-Landau equation (see, e.g., \cite{Guo:1994}) with the following boundary and initial conditions:
\begin{subequations}\label{Landau 2}
\begin{align}
&u_t- a u_{xx}+b u_x+c_1 u+c_2 u^3+c_3 u^5=f(x,t),\ \ (x,t)\in  (0,1)\times \mathbb{R}_+,\\
      &u(0,t)=d_0(t),          u_x(1,t)=d_1(t), \ \ t\in \mathbb{R}_+,\label{Landau 2b}\\
      &u(x,0)=\phi(x),\ \ x\in  (0,1),\end{align}
\end{subequations}
      where 
      {}{$a,c_2,c_3\in\mathbb{R}_+$} and $b,c_1\in\mathbb{R}$ with $ \frac{b^2}{4a}+c_1>0$.

 For \eqref{Landau 2}, we have
 \begin{align*}
&c=c_1,\widetilde{c}=\frac{b^2}{4a}+c_1>0,\\
&h(s)=c_2 s^3+c_3s^5,h'(s)=3c_2 s^2+5c_3s^4\geq 0,\ \forall s\in \mathbb{R},\\
&\alpha_0=1,\alpha_1=0,\beta_0=0,\beta_1=1 \ (\text{note\ that}\ -u_x(0,t)=-d_0(t)),\\
&\widetilde{\alpha}_0=\alpha_0-\frac{b}{2a}\beta_0=1,\widetilde{\beta}_0=\beta_0=0, \\ &\widetilde{\alpha}_1=\alpha_1+\frac{b}{2a}\beta_1=\frac{b}{2a},\widetilde{\beta}_1=\beta_1=1.
\end{align*}
 If we assume further that
  \begin{align}\label{0322+1}
-\frac{2b}{a}< \frac{b^2}{4a}+c_1,\ \ -\frac{b}{2a}\leq a,
\end{align}
 then \eqref{+4b} and \eqref{condition of h} hold. Therefore, \eqref{Landau 2} is EISS. According to \eqref{0320+3} and \eqref{0319+5d} (and noting that $\widetilde{k}_0=0,\widetilde{k}_1>0$ such that $\widetilde{k}_1+\frac{b}{2a}>0$ in \eqref{$k_i$}), it follows {}{that}
 \begin{align*}
&\widetilde{C}_0=0,\widetilde{C}_1=\frac{\widetilde{k}_1\varepsilon_1}{2}-\frac{b}{2a},\\
&\lambda=\widetilde{\lambda}=\min\left\{ \frac{\widetilde{c}}{2}- \frac{\varepsilon}{2} , \frac{\widetilde{c}}{2}-2\widetilde{C}_1- \frac{\varepsilon}{2} \right\},
\end{align*}
 where $ \varepsilon,\varepsilon_1>0$ are small enough such that $\lambda>0$.

  By \eqref{0320+1} and $\widetilde{k}_0=0$, we have
\begin{align*}
\widetilde{\Gamma}(s)=&\frac{4a}{b^2+4ac_1}
\sqrt{\frac{\widetilde{k}_1}{\widetilde{\lambda}\varepsilon_1
}}s+\sqrt{\frac{1 }{\widetilde{\lambda} \varepsilon}}e^{\frac{|b|}{2a}}h\bigg(\frac{4a}{b^2+4ac_1}e^{\frac{|b|}{2a}} s\bigg),\ \forall s\geq 0,\\
\widetilde{\Gamma}_0(s)=&\frac{2a}{b+2a\widetilde{k}_1}
\sqrt{\frac{\widetilde{k}_1}{\widetilde{\lambda}\varepsilon_1
}}
s+
\sqrt{\frac{1 }{\widetilde{\lambda} \varepsilon}}e^{\frac{|b|}{2a}}h\bigg(\frac{2a}{b+2a\widetilde{k}_1}e^{\frac{|b|}{2a}} s\bigg),\ \forall s\geq 0,\\
\widetilde{\Gamma}_1(s)=&\sqrt{\frac{\widetilde{k}_1}{\widetilde{\lambda}\varepsilon_1
}}s+
\sqrt{\frac{1 }{\widetilde{\lambda} \varepsilon}}e^{\frac{|b|}{2a}}h\bigg(e^{\frac{|b|}{2a}} s\bigg),\ \forall s\geq 0,
\end{align*}
which shows that the $\mathcal {K}$-functions in the ISS estimate of \eqref{Landau 2} are given by:
\begin{align*}
&\gamma(s)=\frac{4a}{b^2+4ac_1} e^{\frac{|b|}{2a}}\left( e^{\frac{|b|}{2a}}+
\sqrt{\frac{\widetilde{k}_1}{\widetilde{\lambda}\varepsilon_1
}}\right)s+\sqrt{\frac{1 }{\widetilde{\lambda} \varepsilon}}e^{\frac{|b|}{a}}h\left(\frac{4a}{b^2+4ac_1}e^{\frac{|b|}{2a}} s\right),\ \forall s\geq 0,\\
&{\gamma}_0(s)=e^{\frac{|b|}{2a}}\left( 1+ \frac{2a}{b+2a\widetilde{k}_1}
\sqrt{\frac{\widetilde{k}_1}{\widetilde{\lambda}\varepsilon_1
}}\right)s+\sqrt{\frac{1 }{\widetilde{\lambda} \varepsilon}}e^{\frac{|b|}{a}}h\left(\frac{2a}{b+2a\widetilde{k}_1}e^{\frac{|b|}{2a}} s\right),\ \forall s\geq 0,\\
&{\gamma}_1(s)=e^{\frac{|b|}{2a}}\left(   \frac{2a}{b+2a\widetilde{k}_1} e^{\frac{|b|}{2a}}+ \sqrt{\frac{\widetilde{k}_1}{\widetilde{\lambda}\varepsilon_1
}}\right)s+ \sqrt{\frac{1 }{\widetilde{\lambda} \varepsilon}}e^{\frac{|b|}{a}}h\left(e^{\frac{|b|}{2a}} s\right),\ \forall s\geq 0.
\end{align*}
\begin{remark}{}{It should be mentioned that due to the nonlinear terms in \eqref{Landau 2} and the use of splitting, the gains obtained above are nonlinear and depend on $h$. Moreover, they tend to infinity as $\varepsilon$ (or $\varepsilon_1$) tends to $ 0^+ $. Therefore, it is still a question on how to obtain ISS estimates in $L^2$-norm with uniformly bounded gains for nonlinear PDEs with Neumann boundary disturbances by the Lyapunov method.}
\end{remark}

%
%
%
%


\section{Concluding Remarks}\label{Sec: Conclusion}
This paper presented a new method for the establishment of ISS properties w.r.t. in-domain and boundary disturbances for certain nonlinear parabolic PDEs with different {}{boundary} conditions. The proposed approach for achieving the ISS estimates of the solution is based on the technique of splitting and the weak maximum principle for parabolic PDEs combining with the Lyapunov method. The results show that this method is a convenient tool for the study of ISS properties of PDEs with different disturbances, and it can be applied to stability and regularity analysis for a wider class of nonlinear PDEs with boundary disturbances.

\end{document}